\documentclass[11pt]{article}
\usepackage{mathrsfs}
\usepackage{latexsym,lineno}
\usepackage{epsfig}
\usepackage{color}
\usepackage{amsmath}\usepackage{verbatim}\usepackage{epsf}
\usepackage{amsthm}\usepackage{graphicx, float}\usepackage{graphicx}
\usepackage{amsfonts}\usepackage{amssymb}\usepackage{graphpap}
\usepackage{epic}\usepackage{curves}
\usepackage{tikz}
\usepackage{subfigure}

\newtheorem{fact}{Theorem}

\newtheorem{theorem}[fact]{Theorem}

\newtheorem{definition}{Definition}

\newtheorem{lemma}{Lemma}

\title{\bf Ramsey numbers of graphs \\ with most degrees bounded in random graphs\thanks{Supported in part by NSFC.}}
\author{Ye Wang$\,^{a}$\,\,, \hspace{2mm} Yusheng Li$\,^b$
\footnote{Email:\,ywang@hrbeu.edu.cn,\;li\_yusheng@tongji.edu.cn}\\
{\small $^a\,$College of Mathematical Sciences, Harbin Engineering University, Harbin 150001, China }\\
{\small $^b\,$School of Mathematical Sciences, Tongji University,\; Shanghai 200092, China}  }
\date{}
\begin{document}
\maketitle
\begin{abstract}
For graphs $F$ and $G$, let $F\to G$ signify that any red/blue edge coloring of $F$ contains a monochromatic $G$.
Denote by ${\cal G}(N,p)$ the random graph space of order $N$ and edge probability $p$.
Using the regularity method, one can show that for any fixed $p\in (0,1]$, almost all graphs $F\in{\cal G}(cn,p)$ have
$F\to G$ for any graph $G$ of order $n$ and all but at most $m$ degrees bounded, where $c$ is an integer depending on $p$ and $m$.
Note that $r(K_{m,n})\sim 2^m n$ and $r(K_m+\overline{K}_n)\sim 2^m n$ as $n\to\infty$,
for which we investigate the relation between $c$ and $p$. Let $N=\lfloor c\,2^m n\rfloor$ with $c>1$
and $p_u,p_\ell=\frac{1}{c^{1/m}}(1\pm \sqrt{\frac{M\log n}{n}}\,)$, where $M=M(c,m)>0$.
It is shown that $p_u$ and $p_\ell$ are Ramsey thresholds of $K_{m,n}$ in ${\cal G}(N,p)$.
Namely, almost all $F\in{\cal G}(N,p_u)$ and almost no $F\in{\cal G}(N,p_\ell)$ have $F\to K_{m,n}$, respectively.
Moreover, it is shown that $p_u$ and $p_\ell$ are (ordinary) upper threshold and lower threshold of $K_m+\overline{K}_n$ to appear in ${\cal G}(N,p/2)$, respectively.
We show that ${\cal G}(N,p/2)$ can be identified as the set of red (or blue) graphs obtained from $F\in{\cal G}(N,p)$ by red/blue edge coloring of $F$ with probability $1/2$ for each color,
which leads to the definition of the weak Ramsey thresholds.

\medskip


\medskip

\end{abstract}

\section{Introduction}

~ For graphs $F$ and $G$, let $F\to G$ signify that any red/blue edge coloring of $F$ contains a monochromatic $G$.
The Ramsey number $r(G)$ is defined as the minimum $N$ such that $K_N\to G$.

For two positive functions $f(t)$ and $g(t)$,
let us write $f(t)\le O(g(t))$ or $g(t)\ge \Omega(f(t))$ if  $f(t)\le cg(t)$ for large $t$, where $c>0$ is a constant, and $f(t)=\Theta(g(t))$ if $\Omega(g(t))\le f(t)\le O(g(t))$.
As usual, $f(t)=o(1)$  and $f(t)\sim g(t)$ signify $f(t)\to 0$ and $\frac{f(t)}{g(t)}\to 1$ as $t\to\infty$, respectively.

Let $\Delta(G)$ and $\delta(G)$ be the maximum degree and the minimum degree of a graph $G$, respectively.
The following result follows from a remarkable application of Szemer\'edi's regularity lemma
by Chvat\'al, R\"odl, Szemer\'edi and Trotter \cite{Chvatal}.

\begin{theorem}[\cite{Chvatal}]\label{rl}
Let $\Delta\ge 1$ be a fixed integer. Then there is a constant $c=c(\Delta)>0$ such that  $r(G)\le c n$ for any graph $G$ of order $n$ with $\Delta(G)\le \Delta$.
\end{theorem}
Note that $K_N\to G$ is equivalent to $r(G)\le N$.
The method for the proof is now called the regularity method,
which can be easily generalized to have a result for sparse graphs $F$ as $F\to G$ instead of $K_N\to G$.
It is more convenient to consider sparse graphs in random graphs. Let $V$ be a labelled vertex set with $|V|=N$.
The classic random graph space ${\cal G}(N,p)$ was defined by Erd\H{o}s and R\'enyi \cite{erd-ren},
in which each edge of $K_N$  has a probability $p$ to appear, randomly and independently.
For a graph property $P$, if $\Pr[F\in{\cal G}(N,p):\,F\;\mbox{has}\;P]\to 1$ as $N\to \infty$,
then we say that almost all graphs in ${\cal G}(N,p)$ have the property $P$.

Denote by ${\cal F}_{\Delta,m}$ the family of graphs $F$, where all but at most $m$ vertices of $F$ have degrees at most $\Delta$.

\begin{theorem}\label{t1}
Let $\Delta$ and $m$ be positive integers and $p\in (0,1]$. Then there exists a positive integer $c=c(\Delta,m,p)$
such that almost all graphs $F\in {\cal G}(cn,p)$ have $F\to G$, where $G$ is a graph of order $n$ in ${\cal F}_{\Delta,m}$.
\end{theorem}

It is easy to see the constant $c$ is large if $p$ is small.
In fact, the constant $c$ in the proof of Theorem \ref{t1} is huge.
We shall investigate the relation between  $c$ and $p$.
Since ${\cal F}_{\Delta,m}$ is a large family that contains graphs of many types,
we shall concentrate on some sub-families.

Let $K_{m,n}$ be the complete bipartite graph on $m$ and $n$ vertices and $K_m+\overline{K }_n$ the graph obtained from $K_m$
by adding new edges to connect $K_m$ and other additional $n$ vertices completely.
It is known that $r(K_{m,n})\sim 2^mn$ as $n\to\infty$, but $r(K_m+\overline{K}_n)$ is much harder to handle.
Recently, Conlon \cite{conlon} proved
\[
r(K_m+\overline{K}_n)=(2^m+o(1)) n,
\]
for which the small term $o(1)$ is shown to be at most $\Theta[(\log\log\log n)^{-1/25}]$ in Conlon, Fox and Wigderson \cite{cfw}.
The graphs $K_{m,n}$ and $K_m+\overline{K}_n$ are well-known graphs in ${\cal F}_{\Delta,m}$.

The Ramsey upper threshold $p_u=p_u(N)$ and Ramsey lower threshold $p_\ell=p_{\ell}(N)$ of $G$ in ${\cal G}(N,p)$ are defined by
\[
\lim_{N\to \infty} \Pr[F\in{\cal G}(N,p):\,F\to G]= \left\{  \begin{array}{cl}
1  & \mbox{if $p\ge p_u$,} \\  0 & \mbox{if $p\le p_{_\ell}$}. \end{array}      \right.
\]
For $N=c\,r(K_{m,n})$ with $c>1$, the following result reveals the relation between $c$ and $p$ precisely.
\begin{theorem}\label{rt}
Let $m\ge 1$ be an integer and $N=\lfloor c\,2^m n\rfloor$ with $c>1$.
Then $p_u=\frac{1}{c^{1/m}}\big(1+\frac{\omega(n)}{n}\,\big)$ and $p_\ell=\frac{1}{c^{1/m}}\big(1-\sqrt{{\frac{{M\log n}}{n}}}\,\big)$
are a Ramsey upper threshold and a Ramsey lower threshold of $K_{m,n}$ in ${\cal G}(N,p)$, respectively, where $\omega(n)\to\infty$ and $M=M(c,m)>0$ is a constant.
\end{theorem}
Failing to find a reasonable Ramsey upper threshold of $K_m+\overline{K}_n$,
we shall weaken the property $F\to G$ by generalizing the deterministic event $F\to G$ to a random event,
for which we shall clarify notation.

Let $V$ be a labelled vertex set, and $F$ a fixed graph on $V$. Consider a red/blue edge coloring of $F$ with probability $1/2$ for each color, randomly and independently.
Let $F_R^{1/2}$ and $F_B^{1/2}$ be the red and blue graphs of $F$ in such a coloring,
which are random graphs as the coloring is random even though $F$ is fixed.
Note that different graphs $F$ and $H$ on the same $V$ may yield the same monochromatic graph as $F^{1/2}_R=H^{1/2}_R$ or $F^{1/2}_B=H^{1/2}_B$ with some positive probability.
For a fixed graph $L$ on $V$, if $E(L)\subseteq E(F)$, then there is a unique edge coloring of $F$ such that $F_R^{1/2}=L$
that is a random event.
If $F\in{\cal G}(N,p)$, then $F_R^{1/2}=L$ is a composite event.
Intuitively, each edge of $F_R^{1/2}$ for $F\in{\cal G}(N,p)$ appears with probability $p/2$, and we shall show that this is true.

Recall that for a fixed graph $L$ on labelled vertex set $V$, the probability
\[
\Pr[F\in{\cal G}(N,p):F=L]=p^{e(L)}(1-p)^{\binom{N}{2}-e(L)},
\]
which is said to be the probability of $L$ to appear in ${\cal G}(N,p)$.

\begin{lemma}\label{new-1}
Let $V$ be a labelled vertex set with $|V|=N$ and $L$ a graph on $V$. Then
\[
\Pr[F\in{\cal G}(N,p):F_R^{1/2}= L]= \big(\frac{p}{2}\big)^{e(L)}\big(1-\frac{p}{2}\big)^{\binom{N}{2}-e(L)}.
\]
Namely, the red graphs $F_R^{1/2}$ form ${\cal G}(N,p/2)$ for $F\in{\cal G}(N,p)$.
\end{lemma}

The same argument can be applied for blue graphs $F_B^{1/2}$ as they also form ${\cal G}(N,p/2)$.

For a graph $G$ whose vertex set is not necessarily $V$, we define $F\stackrel{1/2}{\longrightarrow}G$ for $F\in{\cal G}(N,p)$ an event as
\[
\{F\in{\cal G}(N,p):F\stackrel{1/2}{\longrightarrow}G\} = \{F\in{\cal G}(N,p):G\subseteq F_R^{1/2}\;\mbox{or}\;G\subseteq F_B^{1/2}\},
\]
and thus
\begin{eqnarray*}
\{F\in{\cal G}(N,p):G\subseteq F_R^{1/2}\} &\subseteq&\{F\in{\cal G}(N,p):F\stackrel{1/2}{\longrightarrow}G\} \\ &\subseteq&\{F\in{\cal G}(N,p):G\subseteq F_R^{1/2}\} \cup \{F\in{\cal G}(N,p):G\subseteq F_B^{1/2}\}.
\end{eqnarray*}
Note that if $G\subseteq F_R^{1/2}$ and $G\subseteq F_B^{1/2}$, then $F$ are contained in both of the last two events.
In the proof of Lemma \ref{new-1}, we have $\Pr[F_R^{1/2}=L]=\Pr[F_B^{1/2}=L]$ for any fixed graph $L$ on $V$, hence
\[
\Pr\big[F\in{\cal G}(N,p):G\subseteq F_R^{1/2}\big]=\Pr\big[F\in{\cal G}(N,p):G\subseteq F_B^{1/2}\big],
\]
and thus
\[
\Pr[F\in{\cal G}(N,p/2):G\subseteq F]\le\Pr[F\in{\cal G}(N,p):F\stackrel{1/2}{\longrightarrow}G] \le 2\Pr[F\in{\cal G}(N,p/2):G\subseteq F],
\]
which leads to the following definition.

\begin{definition}\label{weak}
For a graph $G$, we call probability functions $p_u=p_u(N)$ and $p_\ell=p_{\ell}(N)$ a weak Ramsey upper threshold and a weak Ramsey lower threshold of $G$ in ${\cal G}(N,p)$, respectively,  if
\[
\lim_{N\to \infty} \Pr[F\in{\cal G}(N,p/2):\,G\subseteq F]= \left\{  \begin{array}{cl}
1  & \mbox{for\, $p\ge p_u$,}\\
0 & \mbox{for\, $p\le p_{_\ell}$}.
\end{array}      \right.
\]
\end{definition}

Even though $F\stackrel{1/2}{\longrightarrow}G$ for $F\in{\cal G}(N,p)$ is a random event,
we do not call the new thresholds as ``random Ramsey thresholds" to avoid confusion since ``thresholds" have their own random meanings,
and we call them ``weak Ramsey thresholds" by viewing  $F\stackrel{1/2}{\longrightarrow}G$ as a weak form of $F\to G$.
In following result, the weak Ramsey thresholds of $K_{m,n}$ are exactly its Ramsey thresholds in Theorem \ref{rt}.
It is interesting to find such relation for general graphs.

\begin{theorem}\label{rt1}
Let $m\ge 1$ be an integer and $N=\lfloor c2^m n\rfloor$ with $c>1$.
Then $\frac{1}{c^{1/m}}\big(1+\frac{\omega(n)}{n}\,\big)$ and $\frac{1}{c^{1/m}}\big(1-\sqrt{{\frac{{M\log n}}{n}}}\,\big)$
are a weak Ramsey upper threshold and a weak Ramsey lower threshold of $K_{m,n}$ in ${\cal G}(N,p)$, respectively, where $\omega(n)\to\infty$ and $M=M(c,m)>0$ is a constant.
\end{theorem}

We give weak Ramsey thresholds of $K_m+\overline{K}_n$ instead of Ramsey thresholds as follows.
\begin{theorem}\label{thr}
Let $m\ge 1$ be an integer and $N=\lfloor c2^m n\rfloor$ with $c>1$.
Then $\frac{1}{c^{1/m}}\big(1\pm\sqrt{\frac{M\log n}{n}}\,\big)$ are a weak Ramsey upper threshold and a weak Ramsey lower threshold
of $K_m+\overline{K}_n$ in ${\cal G}(N,p)$, respectively, where $M=M(c,m)>0$ is a constant.
\end{theorem}

Note that Theorem \ref{rt} implies $r(K_{m,n})\le (2^m+o(1))n$ as $F\to K_{m,n}$ implies $K_N\to K_{m,n}$ for $F$ of order $N=\lfloor c2^m n\rfloor$
and any $c>1$.
However, Theorem \ref{thr} does not give the similar upper bound for $r(K_m+\overline{K}_n)$ as it implies only
there are monochromatic $K_m+\overline{K}_n$ in almost all red/blue edge colorings of $K_N$ instead of any such coloring.

\section{Outline of the proof of Theorem \ref{t1}}

~~As the proof of Theorem \ref{t1} is similar to the proof of Theorem \ref{rl} in \cite{Chvatal},
so we just outline it.

The edge density of a graph $F$ of order $N$ is defined as $d(F)=e(F)/{N\choose 2}$,
in which $e(F)$ denotes the number of edges of $F$.
For a disjoint pair $(X,Y)$ of nonempty disjoint vertex subsets $X$ and $Y$ of $V(F)$,
let $e_F(X,Y)$ be the number of edges of $F$ between $X$ and $Y$.
The ratio $d_F(X,Y)=\frac{e_F(X,Y)}{|X||Y|}$ is called the edge density of $F$ in $(X,Y)$,
which is the probability that any pair $(x,y)$ selected randomly from $X\times Y$ is an edge.
If no danger of confusion, we write $e(X,Y)$ and $d(X,Y)$ for $e_F(X,Y)$ and $d_F(X,Y)$, respectively.

There are many interesting properties of $\mathcal{G}(N,p)$,
see Alon and Spencer \cite{alo-spe}, Bollob\'as \cite{Bollobas}, and Janson, {\L}uczak and Ruci\'nski \cite{jlr}.
Thomason \cite{Thomason1} introduced the notation of jumbled graphs,
and Chung, Graham and Wilson \cite{Chung} showed that many properties are equivalent,
which are satisfied by almost all graphs in ${\cal G}(N,p)$. Some equivalent properties imply the existence of the following graphs.
\begin{lemma}\label{qrp}
Let $p\in(0,1]$ be fixed. Then almost all graphs $F$ in ${\cal G}(N,p)$ have
$d(X,Y)\ge \frac{p}{2}$, where $X$ and $Y$ are disjoint vertex sets and $|X||Y|=\Theta(N^2)$ as $N\to\infty$.
\end{lemma}

Let $\epsilon>0$ be a real number.
A disjoint pair $(X,Y)$ is called $\epsilon$-regular if any $X'\subseteq X$ and $Y'\subseteq Y$
with $|X'|>\epsilon |X|$ and $|Y'|>\epsilon |Y|$ have $|d(X,Y)-d(X',Y')|\le \epsilon$.

The first form of the regularity lemma introduced by Szemer\'edi was in \cite{Szemeredi},
and its general form was in \cite{Szemeredi1}. We need the following formulations of the regularity lemma.

\begin{lemma}[Multi-color regularity lemma]\label{l2}
For any real $\epsilon>0$ and any positive integers $\ell$ and $s$,
there exists $M=M(\epsilon,\ell,s)>\ell$ such that if the edges of a graph $F$ on $N\ge \ell$ vertices are colored in $s$ colors,
all monochromatic graphs have a same partition
$C_1,C_2,\ldots,C_k$ with $\big||C_i|-|C_j|\big|\le 1$, and $\ell\le k\le M$
such that all but at most $\epsilon k^2$ pairs $(C_i,C_j)$, $1\le i<j\le k$,
are $\epsilon$-regular in each monochromatic graph.
\end{lemma}

We shall use Lemma \ref{qrp} and multi-color regularity lemma for the proof of Theorem \ref{t1}
that is similar to the proof of Theorem \ref{rl} in \cite{Chvatal} and thus we omit it.

\section{Proofs of Theorem \ref{rt} -- Theorem \ref{thr}}

~ Let us prove Lemma \ref{new-1} first.
\medskip

{\bf Proof of Lemma \ref{new-1}.} For fixed graphs $F$ and $L$ on the labelled vertex set $V$, if $L\subseteq F$ or equivalently $E(L)\subseteq E(F)$, then
\[
\Pr[F_R^{1/2}=L]=\big(\frac{1}{2}\big)^{e(L)}\big(1-\frac{1}{2}\big)^{e(F)-e(L)}= \frac{1}{2^{e(F)}},
\]
and otherwise $\Pr[F_R^{1/2}=L]=0$. Hence, for a fixed graph $L$ on $V$ and $0<p<1$, we have
\begin{eqnarray*}
\Pr[F\in{\cal G}(N,p):\,F_R^{1/2}=L] &= &  \sum_{F:\;L\subseteq F} \Pr[F\in{\cal G}(N,p)]\Pr[F_R^{1/2}=L] \\
&= & \sum_{F:\;L\subseteq F} p^{e(F)}(1-p)^{\binom{N}{2}-e(F)} \frac{1}{2^{e(F)}} \\
& = & (1-p)^{\binom{N}{2}}\Big[\frac{p}{2(1-p)}\Big]^{e(L)} \sum_{F:\,L\subseteq F} \Big[\frac{p}{2(1-p)}\Big]^{e(F)-e(L)} \\
& = & (1-p)^{\binom{N}{2}}\Big[\frac{p}{2(1-p)}\Big]^{e(L)} \Big[1+\frac{p}{2(1-p)}\Big]^{\binom{N}{2}-e(L)} \\
&= & \Big(\frac{p}{2}\Big)^{e(L)} \Big(1-\frac{p}{2}\Big)^{\binom{N}{2}-e(L)}
\end{eqnarray*}
as required.  The statement for $p=0$ or $p=1$ is easy to verify, and the proof is completed. \hfill $\square$

\medskip

The rest of this section is devoted to the proofs of Theorem \ref{rt} and Theorem \ref{thr},
for which we shall employ the following Chernoff bound \cite{cher,alo-spe,Bollobas}.

\begin{lemma}\label{thr1}
Let $X_1,X_2,\dots$ be mutually independent random variables as
\[
\Pr(X_i=1)=p, \;\; \Pr(X_i=0)=q,
\]
 for any $i$ and $S_T=\Sigma_{i=1}^T X_i$.
Then there exists $\delta_0=\delta_0(p)>0$ such that if $0\le \delta\le \delta_0$,
it holds
\[
\Pr(S_T\ge T(p+\delta))\le\exp\{-T\delta^2/(3pq)\},
\]
and
\[
\Pr(S_T\le T(p-\delta))\le\exp\{-T\delta^2/(3pq)\}.
\]
\end{lemma}

\medskip

{\bf Remark.} In following content, we shall  write $N=\lfloor c2^m n \rfloor$ for fixed integer $m\ge 1$ and constant $c>1$,
and let $\omega(n)$ be a function with $\omega(n)\to\infty$ slowly such that the associated function $p_u\in[0,1]$.
To simplify the notation in the proofs,  we shall admit $c2^m n$ as an integer hence $N=c2^m n$.

\begin{lemma}\label{u-r}
Let $p_u=\frac{1}{c^{1/m}}\big(1+\frac{\omega(n)}{n}\;\big)$.
Then $p_u$ is a Ramsey upper threshold of $K_{m,n}$ in ${\cal G}(N,p)$.
\end{lemma}
{\bf Proof.} ~Let $p_0= \frac{1}{c^{1/m}}\big(1+\frac{\omega(n)}{2n}\;\big)$. We {\bf claim} that almost all graphs $F$
in ${\cal G}(N,p_u)$ have $e(F)> p_0\binom{N}{2}$.

~~Let $T=\binom{N}{2}$ and label all edges of $K_N$ as $e_1,e_2,\dots,e_T$.
For $i=1,2,\dots,T$, define a random variable $X_i$ such that $X_i=1$
if $e_i$ is an edge of $F\in{\cal G}(N,p_u)$ and $0$ otherwise.
Let $e(F)=\sum _{i=1}^T X_i$ be the number of edges of $F$. Note the events that edges appear are mutually independent,
so $e(F)$ has the binomial distribution $B(T,p_u)$.

Let $\delta=p_u-p_0 =\frac{\omega(n)}{2c^{1/m}n}$. By the Chernoff bound,
\begin{eqnarray*}
\Pr[F\in{\cal G}(N,p_u):e(F)\le p_0 T] & = & \Pr[F\in{\cal G}(N,p_u):e(F)\le T (p_u-\delta)] \\
& \le &\exp\{-T\delta^2/(3p_u(1-p_u))\},
\end{eqnarray*}
in which $3p_u(1-p_u)\to \frac{3(c^{1/m}-1)}{c^{2/m}}>0$ and
\[
T\delta^2\sim\frac{c^2 2^{2m} n^2}{2}\cdot\frac{\omega^2(n)}{4c^{2/m}n^2}=2^{2m-3}c^{2(1-1/m)}\omega^2(n) \to \infty.
\]
So the evaluated probability tends to zero, which implies that almost all graphs $F\in{\cal G}(N,p_u)$ have $e(F)> p_0\binom{N}{2}$,
and the claim follows.

Then, we shall show that
\begin{equation}\label{up-ram}
\Pr[F\in{\cal G}(N,p_u):\,F\to K_{m,n}] \to 1.
\end{equation}
From double counting of K\"ovari, S\'os and Tur\'an \cite{kst}, we know
\begin{equation*}\label{bi-turan}
ex(N;K_{m,n})\le \frac{1}{2}\big[(n-1)^{1/m}N^{2-1/m} + (m-1)N\big].
\end{equation*}
It suffices that $2ex(N;K_{m,n}) < e(F)$ for almost all graphs $F$ in ${\cal G}(N,p_u)$.
Therefore, by the claim that has been shown, we shall show
\[
(n-1)^{1/m}N^{2-1/m} + (m-1)N < \frac{1}{c^{1/m}}\Big(1+\frac{\omega(n)}{2n}\Big)\binom{N}{2}.
\]
It suffices to show that
\[
\Big(\frac{N}{c2^m}\Big)^{1/m} N^{2-1/m} + m N \le \frac{1}{2c^{1/m}} \Big(1+\frac{\omega(n)}{2n}\Big) N(N-1),
\]
equivalently,
\[
N^2 + 2c^{1/m}m N \le  \Big(1+\frac{\omega(n)}{2n}\Big)(N^2 - N),
\]
which follows by deleting square terms on both sides and noting $\omega(n)\to\infty$.  \hfill $\square$

\medskip

We will use a simple fact as follows. For  graph properties $P$ and $Q$,
if $\Pr(P)\to 1$ and $\Pr(Q)\to 1$ in the same random graph space ${\cal G}(N,p)$, then $\Pr(P\cap Q)\to 1$
since $\Pr(P\cap Q)\ge 1- \Pr(\overline{P})-\Pr(\overline{Q})$,
where $\Pr(P)=\Pr[F\in{\cal G}(N,p):F\,\mbox{has}\,P]$.

\begin{lemma}\label{l-r}
Let $p_{_\ell}= \frac{1}{2c^{1/m}}\big(1-\sqrt{\frac{M\log n}{n}}\,\big)$, where $M=M(c,m)>0$ is a constant.
Then $p_{_\ell}$ is a lower threshold of $K_{m,n}$ in ${\cal G}(N,p)$, and $2p_{_\ell}$ is a Ramsey lower threshold of $K_{m,n}$ in ${\cal G}(N,p)$.
\end{lemma}
{\bf Proof.} Note $(1-x)^{1/m}=1-\frac{1}{m}x+o(x^2)$ for $x\to 0$,
hence $1-\frac{2}{m}x\le (1-x)^{1/m}\le 1-\frac{1}{2m}x$ for small $x>0$.
By slightly shifting $M$, it suffices to show that
\[
p_{_\ell}= \frac{1}{2}\Big[\frac{1}{c}\big(1-\sqrt{\frac{M\log n}{n}}\;\big)\Big]^{1/m}
\]
is such a lower threshold, i.e., $\Pr[F\in{\cal G}(N,p_\ell):\,K_{m,n}\subseteq F] \to 0.$

Let $U$ be a fixed subset of $V$ with $|U|=m$. Let $v_1,v_2,\dots,v_{N-m}$ be the vertices outside of $U$.
For each $i=1,2,\dots, N-m$, define a random variable $X_i$ such that $X_i=1$
if $v_i$ is a common neighbor of $U$ in $F$ and $0$ otherwise. Then $\Pr(X_i=1)=p_{_\ell}^m$.

Set $S_{N-m}=\sum _{i=1}^{N-m}X_i$ that has the binomial distribution $B(N-m,p_0)$, where $p_0=p_\ell^m$.
Note that the event $S_{N-m}\ge n$ means that there is a $K_{m,\,n}$ with $U$ as the part of $m$ vertices.
Hence
\[
\Pr\big[K_{m,n}\subseteq {\cal G}(N,p_{\ell})\big] \le \binom{N}{m}\Pr(S_{N-m}\ge n).
\]
We now evaluate the probability $\Pr(S_{N-m}\ge n)$. Write
\[
n=\frac{N}{c2^m} = (p_0 +\delta)(N-m),
\]
where
\[
\delta= \frac{N}{c2^m(N-m)}-p_0 = \frac{1}{c2^m} \Big[\frac{N}{N-m}-1 +\sqrt{\frac{M\log n}{n}}\Big]=\frac{1}{c2^m}\Big[\frac{m}{N-m}+\sqrt{\frac{M\log n}{n}}\;\Big].
\]
Note $\frac{m}{N-m}\ll \sqrt{\frac{M\log n}{n}}$. By virtue of Chernoff bound,
\[
\Pr(S_{N-m}\ge n) = \Pr\big[S_{N-m}\ge (p_0 +\delta)(N-m)\big]\le \exp\big\{-(N-m)\delta^2/(3p_0q_0)\big\}.
\]
Note that $3p_0q_0$ tends a positive constant, and
\[
(N-m)\delta^2 \sim N \delta^2 \sim c2^m n \frac{M\log n}{c^2 2^{2m} n} = \frac{M}{c 2^m}\log n.
\]
Taking large $M$ such that $(N-m)\delta^2/(3p_0q_0) >2m\log n$ for large $n$, then we have
\[
{N\choose m}\Pr(S_{N-m}\ge n) \le \Theta\Big(\frac{n^m}{n^{2m}}\Big) \to  0,
\]
and thus $p_\ell$ is a lower threshold of $K_{m,n}$ in ${\cal G}(N,p)$.

Now we shall show $2p_{\ell}$  is a Ramsey lower threshold of $K_{m,n}$ in ${\cal G}(N,p)$.
Recall the proof of the claim in Lemma \ref{u-r},
we are easy to see that almost all graphs $F$ in ${\cal G}(N,p_{\ell})$ have
\begin{equation}\label{narrow}
A=\big(p_{\ell}-\frac{\omega(n)}{2c^{1/m}n}\big)\binom{N}{2} \le  e(F)\le \big(p_{\ell}+\frac{\omega(n)}{2c^{1/m}n}\big)\binom{N}{2}=B.
\end{equation}
Denote by $T=2^{^{\binom{N}{2}}}$. Let us label all graphs on labelled vertex set $V$ as $R_1,R_2,\dots,R_T$
and set $t=t(n)$ with $1< t < T$ such that $R_i$ contains no $K_{m,n}$ for $1\le i\le t$ and
\begin{equation}\label{to-zero}
\sum_{i=1}^t \Pr[R_i:R_i\in{\cal G}(N,p_\ell)]\to 1
\end{equation}
as $n\to\infty$. Furthermore, we assume that the number of edges of each $R_i$ with $1\le i\le t$ is in the interval $[A,B]$ described
in (\ref{narrow}).

We then consider red/blue edge coloring of $F$ from ${\cal G}(N,2p_{\ell})$ with probability $1/2$ for each edge, randomly and independently,
for which all the red graphs $F_R$ form the random graph space ${\cal G}(N,p_{\ell})$ by Lemma \ref{new-1}.
The same argument can be applied for blue graphs $F_B$ for $F\in{\cal G}(N,2p_{\ell})$,
which form the same space ${\cal G}(N,p_{\ell})$.

The similar proof implies that almost all $F\in{\cal G}(N,2p_{\ell})$
have $e(F)\in [2A,2B]$, namely $e(F_R)+e(F_B)$ are in the interval $[2A,2B]$.
Combining this and the fact that both $e(F_R)$ and $e(F_B)$ are in the same interval $[A,B]$,
we can further assume that if $F_R=R_i$ for $1\le i\le t$, then there is some $j$ with $1\le j\le t$ such that $F_B=R_j$,
namely $R_i$ and $R_j$ are the red and blue graphs obtained from edge coloring of $F$, respectively.

Therefore, for almost all $F\in{\cal G}(N,2p_{\ell})$, there is a red/blue edge coloring of $F$
such that there is no monochromatic $K_{m,n}$, completing the proof. \hfill $\square$

\medskip

{\bf Proof of Theorem \ref{rt}.} The proof follows from Lemma \ref{u-r} and Lemma \ref{l-r}.  \hfill $\square$

\medskip
{\bf Proof of Theorem \ref{thr1}.}
The weak Ramsey lower threshold of $K_{m,n}$ is from Lemma \ref{l-r},
and the weak Ramsey upper threshold of $K_{m,n}$ follows from the weak Ramsey upper threshold of
$K_m+\overline{K}_n$ in Theorem \ref{thr}.
\hfill $\square$

\medskip

{\bf Proof of Theorem \ref{thr}.}
Since $p_\ell=\frac{1}{c^{1/m}}\big(1-\sqrt{\frac{M\log n}{n}}\,\big)$ is a weak Ramsey lower threshold of $K_{m,n}$ in ${\cal G}(N,p)$,
it is also a weak Ramsey lower threshold of $K_m+\overline{K}_n$.
Then, we shall show that $\frac{1}{c^{1/m}}\big(1+\sqrt{\frac{M\log n}{n}}\,\big)$ is a weak Ramsey upper threshold
of $K_m+\overline{K}_n$ in ${\cal G}(N,p)$, namely almost all graphs ${\cal G}(N,p_u/2)$ contain $K_m+\overline{K}_n$.

The same argument as that in the proof of Lemma \ref{l-r} tells us that  we may assume that
\[
p_{u}= \Big[\frac{1}{c}\big(1 +\sqrt{\frac{M \log n}{n}}\;\big)\Big]^{1/m}.
\]

Let $U$ be a fixed subset of $V$ with $|U|=m$. Let $v_1,v_2,\dots,v_{N-m}$ be the vertices outside of $U$.
For each $i=1,2,\dots, N-m$, define a random variable $X_i(U)$ such that $X_i(U)=1$
if $v_i$ is common neighbor of $U$ and $0$ otherwise. Then $\Pr[X_i(U)=1]=(p_u/2)^m$.

Set $p_0=(p_u/2)^m$ and $S_{N-m}=\sum _{i=1}^{N-m}X_i(U)$.
The the event $S_{N-m}\le n-1$ means that $U$ has at most $n-1$ common neighbors.
Note that $S_{N-m}$ has the binomial distribution $B(N-m,p_0)$ that is independent to $U$.

We now evaluate the probability $\Pr(S_{N-m}\le n-1)$. Write
\[
n-1=\frac{N-c2^m}{c2^m} =\Big[\big(\frac{p_u}{2}\big)^m -\delta\Big](N-m) = (p_0-\delta)(N-m),
\]
where
\begin{eqnarray*}
\delta & = & \frac{1}{2^m} \Big[p_u^m-\frac{N-c2^m}{c(N-m)}\Big]= \frac{1}{c2^m}\Big[\sqrt{\frac{M\log n}{n}}+\frac{c2^m-m}{N-m}\Big].
\end{eqnarray*}
By virtue of Chernoff bound,
\[
\Pr(S_{N-m}\le n-1) = \Pr\big[S_{N-m}\le (p_0-\delta)(N-m)\big]\le \exp\big\{-(N-m)\delta^2/(3p_0q_0)\big\}.
\]
Note that $p_0q_0$ tends a positive constant, and
\[
-(N-m)\delta^2 < - \frac{1}{2} N\frac{M \log n}{c^2 2^{2m} n}= - \log n\cdot M',
\]
where $M'/(3p_0q_0) >2m$ if $M$ is large.
Let $A_U$ be the event that the number of common neighbors of $U$ in $F\in {\cal G}(N,p_u/2)$ is at most $n-1$, where $|U|=m$.
Then
\[
\Pr(A_U)=\Pr(S_{N-m}\le n-1)<\Big(\frac{1}{n}\Big)^{2m}.
\]
Therefore, we have
\[
\Pr\big[\cup_{U} A_U\big]\le \binom{N}{m}\Pr[S_{N-m}\le n-1]\le N^m \Big(\frac{1}{n}\Big)^{2m}\to 0.
\]
Let $Y=\cap_{U}\overline{A}_U$ be the event that each vertex set $U$ of graph $F\in {\cal G}(N,p_u/2)$  has at least $n$ common neighbors, where $|U|=m$.
The $\Pr(Y)\to 1$.

A classic result of Erd\H{o}s and Reny\'i \cite{erd-ren} tells us that $\frac{\omega(N)}{N^{2/(m-1)}}$
is an upper threshold of $K_m$ to appear in ${\cal G}(N,p)$, where $\omega(N)\to \infty$ slowly arbitrarily.
Since $p_u\sim \frac{1}{c^{1/m}}\gg \frac{\omega(N)}{N^{2/(m-1)}}$, we know that $p_u$ is an upper threshold of $K_m$ to appear in ${\cal G}(N,p/2)$.

Denote by $Z$ the event that there is a clique $U$ in  $F\in {\cal G}(N,p_u/2)$ with $|U|=m$.
Then $\Pr(Z)\to 1$ and thus $\Pr(Y\cap Z)\to 1$. So almost all graphs in ${\cal G}(N,p_u/2)$ contain a $K_m+\overline{K}_n$,
and the proof is completed. \hfill $\square$

\medskip

{\bf Acknowledgement}
Recently, Conlon, Fox and Wigderson \cite{cfw1} have determined off-diagonal book Ramsey number $r(B_{an}^{(m)},B_n^{(m)})$ to be $(a^{1/m}+1)^m n+o_m(n)$ as $n\to \infty$, where $0<a\le 1$ and $B_n^{(m)}=K_m+\overline{K}_n$.

\end{document}